\input amstex
\documentstyle{amsppt}
\input bull-ppt
\keyedby{gabai/amh}

\define\cirD{\overset\circ\to D}

\define\cirW{\overset\circ\to W}

\define\Si{S^2_{\infty}}
\define\lambdaij{\lambda_{ij}}
\define\sigmaij{\sigma_{ij}}

\define \finv{f^{-1}}

\define\hyp{\Bbb H}
\define\BR{\Bbb R}
\define\BB{\Bbb B}

\define\Homeo{\operatorname{ Homeo}}
\define\Out{\operatorname{ Out}}
\define\Isom{\operatorname{ Isom}}


\topmatter
\cvol{31}
\cvolyear{1994}
\cmonth{October}
\cyear{1994}
\cvolno{2}
\cpgs{228-232}

\title On the geometric and topological rigidity of hyperbolic 3-manifolds \endtitle
\author David Gabai\endauthor
\address California Institute of Technology, Pasadena, 
California
 91125\endaddress 
\ml Gabai\@math.caltech.edu \endml 
\date November 16, 1993\enddate
\thanks Partially supported by NSF Grants DMS-8902343,  
DMS-9200584,
and SERC grant GR/H60851\endthanks 
\abstract A homotopy equivalence between a hyperbolic 
3-manifold and a closed
irreducible 3-manifold is homotopic to a homeomorphsim 
provided the hyperbolic 
manifold satisfies a purely geometric condition. There are 
no known examples of
hyperbolic 3-manifolds which do not satisfy this 
condition.\endabstract
\endtopmatter %

\document
One of the central problems of 3-manifold topology is to
determine when a homotopy equivalence between two closed 
orientable
 irreducible 3-manifolds is
homotopic to a homeomorphism.  If one of these manifolds 
is $S^3$, 
then this is Poincar\'e's problem.
The results of  [Re], [Fr], [Ru], [Bo], and [HR] (see also 
[Ol])
 completely solve this problem for
maps between lens spaces.  In particular there exist 
nonhomeomorphic
 but homotopy equivalent
lens spaces (e.g. L(7,1) and L(7,2)), and there exist 
self-homotopy 
equivalences not homotopic to
homeomorphisms (e.g. the self-homotopy equivalence of 
L(8,1) whose 
$\pi_1$-map is multiplication by
3).  By Waldhausen [W] (resp. Scott [S]) a homotopy 
equivalence between 
a closed Haken 3-manifold
(resp.\ a Seifert-fibred space with infinite $\pi_1$) and 
an irreducible
 3-manifold can be homotoped
to a homeomorphism. By Mostow [M] a homotopy equivalence 
between two
 closed hyperbolic 3-manifolds
can be homotoped to a homeomorphism and in fact an 
isometry. However,
 the general case of homotopy
equivalence between a hyperbolic 3-manifold and an 
irreducible 3-manifold
 remains to be
investigated.    These problems and results should be 
contrasted with
 the conjecture [T] that a
closed irreducible orientable 3-manifold is either Haken, 
or Seifert
 fibred with infinite $\pi_1$, or
the quotient of $S^3$ by an orthogonal action, or the 
quotient of $
\hyp^3$ via a cocompact group of
hyperbolic isometries.
\proclaim {Theorem 1 \rm [G2]} Let $N$ be a closed, 
orientable, hyperbolic 
\RM3-manifold containing an embedded hyperbolic tube  of  
radius $(\log3
)/2=.549306\dotsc$ about a closed
geodesic.   Then\RM:
\roster
\item"(i)"  If $f:M\to N$ is a homotopy equivalence where 
$M$ is
an irreducible \RM3-manifold, then $f$
is  homotopic to a homeomorphism.
\item"(ii)" If $f,g:M\to N$ are homotopic  homeomorphisms, 
then $f$ is 
isotopic to $g$. 
\item"(iii)" The space of hyperbolic metrics on $N$ is 
path connected.    
\endroster
\endproclaim
\demo {Remarks}  (i) Thus $N$ is both topologically and 
geometrically 
rigid provided it
satisfies the purely geometric condition of having a 
modest-sized 
tube about a geodesic. 
Actually the conclusion of Theorem 1 holds provided $N$ 
satisfies a more general
geometric/topological property called the {\it insulator 
condition}.

(ii) Jeff Weeks's tube radius/ortholength program \cite{We}
has found one
hyperbolic  3-manifold $N$ (of volume 1.0149...) which 
fails to have
 a $(\log3)/2$ tube. 
Again by Weeks, conclusions (i)--(iii) above are 
applicable to $N$ because
 it  satisfies the insulator
condition.  \footnote"*"{Note added in proof, June 21, 
1994: Several weeks
ago, Nathaniel Thurston, of the Geometry Center at the 
University of
Minnesota, discovered what appear to be five additional 
hyperbolic 3-manifolds
with a shortest geodesic which does not have a $\log(3)/2$ 
tube. There is
strong evidence that these manifolds satisfy the insulator 
condition. Thurston
made use of Robert Riley's POINCAR\'E program as well as 
some ideas of Robert
Meyerhoff and the author.}

(iii) An application of the hyperbolic law of cosines 
shows that if 
the shortest geodesic $\delta$
in $N$ has length  $>1.353$, then  tube radius 
($\delta)>(\log3)/2$. 

(iv) If $N$ has a geodesic $\delta$ of length $<0.0978$, 
then Meyerhoff\,'s
 tube radius formula \cite{Me, \S3}
implies that tube radius $(\delta)>(\log3)/2$.  Recently 
Gehring and Martin [GM1, 2] improved this number to 
$0.19$.   Combined 
with the work of Jorgenson
[Gr], this shows that for any $n>0$ there exist only 
finitely many
  hyperbolic 3-manifolds of volume
$<n$  which can fail to satisfy  the hypothesis of Theorem 
1.

(v)   Farrell and Jones [FJ] showed that if $f:M\to N$ is 
a homotopy 
equivalence between closed manifolds and $N$ is a
hyperbolic manifold of dimension $\ge 5$, then $f$ is 
homotopic to 
a homeomorphism.\enddemo
 
The theme of the proof of Theorem 1 is to abstract the 
ideas in the
 following example to the
setting  of homotopy hyperbolic 3-manifolds.
\ex{Example 2}   Let $\delta$ be a simple closed geodesic 
in the 
hyperbolic 3-manifold $N$.  $\delta$ lifts to a collection 
$\Delta=\{\delta_i\}$ of hyperbolic lines in $\hyp^3$.  To 
each pair 
$\delta_i,\delta_j$ there exists the {\it midplane} 
$D_{ij}$, i.e. the 
hyperbolic halfplane orthogonal to and cutting the middle 
of the {\it 
orthocurve}
(i.e. the shortest line segment) between $\delta_i$ and 
$\delta_j$.  
Each
$D_{ij}$  extends to a circle $\lambda_{ij}$ on 
$S^2_{\infty}$, which
separates  $\partial \delta_i$ from $\partial\delta_j$.  
Now fix $i$.  
Let
$H_{ij}$ be  the closed $\hyp^3$-halfspace bounded by 
$D_{ij}$ 
containing
$\delta_i$.  $W_i=\cap  H_{ij}=D^2\times \BR$  is the {\it 
Dirichlet
 tube domain}
associated to the geodesic  $\delta_i$. $\cirW_i$ projects 
to an open 
solid
torus in $N$ containing $\delta$ as its core.   In fact,
 $W=W_i/\langle \delta_i\rangle$ is a
solid torus with boundary a finite union of  totally 
geodesic polygons.  
\endex
\dfn{Definition 3}  Let $\Cal A=\{\lambdaij\}$.  We call 
the pair
 $(\pi_1(N),\Cal A)$ the 
{\it Dirichlet insulator family} associated to $\delta$.  
It is 
{\it  noncoalesceable}  if for no
$i$ does there exist $\lambda_{ij_1},\lambda_{ij_2},  
\lambda_{ij_3}$
 whose union separates the
points of  $\partial\delta_i$. A hyperbolic 3-manifold 
satisfies the
 {\it insulator condition}
if  the Dirichlet insulator family associated to some 
geodesic is 
noncoalesceable.  \enddfn
\proclaim {Conjecture 4}  The Dirichlet insulator family 
associated
 to a shortest geodesic in a
closed orientable hyperbolic \RM3-manifold is 
noncoalesceable.     \endproclaim
\demo {Remarks} (i)  The notion of Dirichlet insulator can 
be abstracted 
to the more general
notion of {\it insulator}  [G2].  

(ii)  If a geodesic $\delta$ has tube radius $> (\log 
3)/2$, 
then 
its Dirichlet insulator family is noncoalesceable. The 
explanation 
boils down to the following
observation in 2-dimensional hyperbolic geometry.   If a 
geodesic $
\gamma_1$ is at distance
log(3)/2 from a point $x$, then in the visual circle of 
$x$,  
$\gamma_1$ takes up exactly 120
degrees.  Thus three geodesics of distance $>$log(3)/2 
from $x$ cannot
 form a link around  $x$. 
\enddemo
\proclaim {Proposition 5 \rm [G1, 2]} If $f:M\to N$ is a 
homotopy 
equivalence between the closed hyperbolic \RM3-manifold 
$N$ and the 
irreducible \RM3-manifold $M$, then 
$M$ and $N$ are covered by the same closed hyperbolic 
manifold $X$.
The covering map $p_1:X\to M$ can be chosen so that the 
homotopy 
equivalence lifts and
extends to a  mapping $f_1:X\to X$ homotopic to 
$\roman{id}_X$.  The group 
of covering transformations on
$\hyp^3$
 defined by $\pi_1(M)$ and $\pi_1(N)$ induce identical 
group actions
 on $\Si$. \endproclaim
The following proposition gives a criterion for showing 
that a homotopy 
equivalence can be deformed
into a homeomorphism.  
\proclaim {Proposition 6}  Let $f:M\to N$ be a homotopy 
equivalence 
between the closed orientable hyperbolic \RM3-manifold $N$ 
and the 
irreducible \RM3-manifold 
$M$.
If there exists a simple closed curve  $\gamma\subset M$, a 
geodesic $\delta \subset N$ and a homeomorphism 
$k:(\BB^3,p^{-
1}(\gamma))\to 
(\BB^3,q^{-1}(\delta))$ such that $k\mid \partial \BB^3 = 
\roman{id}$, then 
$f$ is 
homotopic to 
a homeomorphism. \endproclaim

\topspace{16.5pc}

\demo {Remarks} (i)  Implicit in the statement of 
Proposition 6 is 
the definition of
$p,q$ and the identification of $\Si$'s given in 
Proposition 5.

(ii)  Said another way, $f$ is homotopic to a
homeomorphism provided the $\BB^3$-link
$\Delta$ is  equivalent to 
the $\BB^3$-link $\Gamma$.
Here $\Delta$ is the preimage of $\delta$ in $\hyp^3$ 
extended to 
$\BB^3$, and $\Gamma$ is defined similarly.  That these 
links are 
equivalent means 
that there exists a homeomorphism $k:(\BB^3,\Gamma) \to 
(\BB^3,\Delta)$ so that $k\mid S^2_{\infty}=$ id. 
\proclaim {Theorem 7}  Let $f:M\to N$ be a homotopy 
equivalence, 
where 
$M$ is a closed, connected, orientable, irreducible 
\RM3-manifold and 
$N$ is a 
hyperbolic \RM3-manifold.  If $N$ 
possesses a geodesic $\delta$ with a noncoalesceable 
insulator family, 
then $f$ is homotopic to a homeomorphism.
\endproclaim
\enddemo
\demo {Outline of the proof}  To each smooth simple closed 
curve 
$\lambda_{ij}$ in  $S^2_{\infty}$, there exists a 
lamination 
$\sigma_{ij}$ by 
least area (with respect to the
metric induced by $M$) planes in $\hyp^3$,  with limit set  
$\lambda_{ij}$ such that $\sigma_{ij}$ lies in a fixed 
width hyperbolic regular 
neighborhood of the hyperbolic convex 
hull of $\lambda_{ij}$.  
Here $\{\lambda_{ij}\}$ is the 
$(\pi_1(N),\{\partial\delta_i\})$ and 
hence $(\pi_1(M),\{\partial\delta_i\})$ noncoalesceable 
insulator 
family.
Fix $i$.  Let
$H_{ij}$ be the $\hyp^3$-complementary region of $\sigmaij$ 
containing the ends of $\delta_i$.
We show that $\cap_j H_{ij}$ contains a component 
$V_i=\cirD^2 \times \BR$
which projects to an open  solid torus in $M$.  Define  
$\gamma$ to
 be the core of this solid torus 
and
$\gamma_i$ the lift which lives in $V_i$.  The isotopy 
class of 
$\gamma$ is independent of all choices, i.e. the metric on 
$M$ and 
the choice of
$\{\sigma_{ij}\}$ for a fixed  metric.  Let  $\tau_0$ be 
the link in $X$ 
which is
the preimage of $\gamma$, so $\{\gamma_i\}$ is also the 
set of lifts
 of components of
$\tau_0$ to
$\hyp^3$.  The Riemannian metric
$\mu_0$ on $X$ induced from $M$ and the hyperbolic metric 
$\mu_1$ on $X$ are
connected by a smooth path  $\mu_t$ of metrics.  These 
metrics lift 
to
$\pi_1(X)$ equivariant metrics $\tilde\mu_t$ on  $\hyp^3$, 
so the 
above
construction  applied to the  
$(\pi_1(X),\{\partial\delta_i\})$ 
insulator family
$\{\lambda_{ij}\}$ with respect
to the  $\tilde\mu_t$ metric  yields   a link $\tau_t$  in 
$X$.  Since 
the
isotopy class of $\tau_t$ is independent of $t$,  $\tau_0$ 
is
isotopic to  $\tau_1$, the preimage of $\delta$ in $X$.  
We conclude 
that the
$\BB^3$-link $\Gamma$ is  equivalent to the $\BB^3$-link 
$\Delta$ ,
and so by Proposition 6
$f$ is  homotopic to a homeomorphism.  \enddemo
\proclaim {Theorem 8}  If $N$ is a closed, oriented, 
hyperbolic 
\RM3-manifold possessing a geodesic
$\delta$  with a noncoalesceable insulator family and  
$f:N\to N$ is a 
homeomorphism homotopic to {\rm id}, then $f$ 
is 
isotopic to {\rm id}.\endproclaim
\demo {Idea of the proof}  Let $\rho_0$ denote the 
hyperbolic 
metric  on $N$.  Let $\rho_1$ be the pull-back 
hyperbolic metric on $N$ induced via $f$, which we can 
assume 
is a diffeomorphism.  These metrics 
are 
connected by a family $\rho_t$.  As in 
the proof of Theorem 7, to each $\rho_t$ there is 
associated a 
simple closed curve $\gamma_t$ where 
$\gamma_0=\delta$ and $\gamma_1=\finv(\delta)$ and all of 
these $\gamma_t$'s are isotopic.  
Therefore, $f$ is isotopic to a map which fixes $\delta$ 
pointwise.  
A theorem of Siebenmann [BS] implies that $f$ is isotopic 
to 
id.  \enddemo
\proclaim {Corollary 9} If $N$ satisfies the insulator 
condition, then 
$$\Homeo(N)/\Homeo_0(N) = \Out(\pi_1(N)) = 
\Isom(N)\.$$\endproclaim
\demo {Proof} Since a hyperbolic 3-manifold is a 
K($\pi$,1), homotopy classes of
homeomorphisms are  parametrized by Out($\pi_1(N)$).  
Mostow implies that each
homotopy class is representable by a  unique 
isometry.  Theorem 8 implies that homotopy classes of 
homeomorphisms
 are the same as isotopy
classes  of homeomorphisms.\qed\enddemo
\demo {Remark \rm 10  (Why a coalesceable insulator is 
bad)} It is possible
that the $\cap_jH_{ij}$  resulting from the construction 
applied to
 a coalesceable
insulator family
 would  contain no  
$\cirD^2\times \BR$ component.  In fact, using the wrong 
metric,  some 
hyperbolic plane $P$ 
transverse to $\delta_i$ may be disjoint from 
$\cap_jH_{ij}$.   
This is the usual
problem of a triangular prism formed by three minimal 
surfaces being 
obliterated upon change of
metric.\enddemo 
\demo {Remark \rm 11 (What Mostow does not say)}  If 
$\rho_0$ is a hyperbolic
metric on $N$, then  a nontrivial element $\alpha$ of 
$\pi_1(N) $ 
determines a geodesic
$\delta_0$ on $N$.  Mostow's rigidity theorem does not 
rule out the
  possibility that with respect 
to a different hyperbolic metric $\rho_1$,  the geodesic 
$\delta_1$
 associated to $\alpha$ would 
lie in a 
different isotopy class than $\delta_0$.  What Mostow does 
assert is that there
exists a diffeomorphism $f:N\to N$, {\it  homotopic} to 
id, such that
$f(\delta_0)=\delta_1$.  Theorems 7 and 8 show that if $N$ 
satisfies
the insulator condition, then $\delta_0$ is isotopic to 
$\delta_1$ 
and, further, that  the 
diffeomorphism is isotopic to id. 
Said another way, Mostow asserts that hyperbolic 
structures are unique up to
homotopy, while Theorem 1 asserts that under mild 
hypothesis a hyperbolic
structure is unique up to isotopy.
\enddemo

\heading Acknowledgments\endheading

Special thanks to Charlie Frohman, Joel 
Hass, Robert Meyerhoff, Peter
Scott,  Evelyn Strauss, and the Mathematics
Institute of the  University of Warwick. 
\Biblio
\widestnumber\key{GM1}
\ref\key Bo\by F. Bonahon\paper Diffeotopes des espaces 
lenticulaires
\jour Topology\vol
22\yr 1983\pages 305--314\endref
\ref\key BS\by F. Bonahon and L. Siebenmann\toappear\endref
\ref\key FJ\by F. T. Farrell and L. Jones\paper A 
topological
analogue of Mostow\RM's rigidity theorem\jour J. Amer. 
Math. Soc. \vol
 2\yr 1989\pages
257--370\endref
\ref\key Fr\by W. Franz\paper Abbildungsklassen und 
fixpunktklassen
 dreidimensionalen linsenraume
\jour J. Reine. Angew. Math.\vol 185\yr 1943\pages 
65--77\endref
\ref \key G1 \by Gabai
\paper Homotopy hyperbolic \RM3-manifolds are virtually 
hyperbolic
\jour  J. Amer. Math. Soc.
\vol 7\yr1994\pages 193--198
\endref
\ref \key G2 \bysame
\paper On the geometric and topological rigidity of 
hyperbolic \RM3-manifolds
\paperinfo preprint
\endref
\ref \key GM1\by F. Gehring and G. Martin\paper 
Commutators, 
collars and the geometry of Mobius groups\jour J. d'Analyse
\toappear \endref
\ref\key GM2\bysame\paper Torsion and volume in 
hyperbolic \RM3-folds\paperinfo in preparation\endref
\ref \key Gr \by M. Gromov \paper Hyperbolic manifolds 
according to Thurston
and Jorgensen \jour  Sem. Bourbaki \vol 32 \yr1979 \pages 
40--52
\endref
\ref \key Me \by R. Meyerhoff
\paper A lower bound for the volume of hyperbolic 
\RM3-manifolds
\jour Canada. J. Math.  \vol 39 \yr1987 \pages 1038--1056
\endref
\ref \key Mo \by G. D. Mostow
\paper Quasiconformal mappings in n-space and the rigidity 
of hyperbolic
 space forms
\jour Inst. Hautes \'Etudes Sci. Publ. Math.
\vol 34 \yr1968 \pages 53--104
\endref
\ref\key Ol\by P. Olum\paper Mappings of manifolds and the 
notion 
of degree\jour Ann. of Math. (2)
\vol 58\yr 1953\pages 458--480\endref
\ref\key Re
\by K. Reidemeister
\paper Homotopieringe und Linsenraume
\jour Abh. Math. Sem. Univ. Hamburg
\vol 11
\yr 1935
\pages 102--109
\endref
\ref\key Ru\by M. Rueff\paper Beitrage zur untersuchung 
der abbildungen
 von mannigfaltigkeiten
\jour Compositio Math.\vol 6\yr 1938\pages 161--202\endref
\ref \key S \by P. Scott
\paper There are no fake Seifert fibred spaces with 
infinite $\pi_1$
\jour Ann. of Math.  (2) \vol 117 \yr1983 \pages 35--70
\endref
\ref \key T \by W. P. Thurston
\paper Three-dimensional manifolds, Kleinian groups, and 
hyperbolic geometry
\inbook Proc. Sympos. Pure Math., vol. 39
\publ Amer. Math. Soc.
\publaddr Providence, RI\yr1983 \pages 87--111\endref
\ref \key W  \by F. Waldhausen
\paper On irreducible \RM3-manifolds which are 
sufficiently large
\jour Ann. of Math. (2) \vol 87 \yr1968 \pages 56--88
\endref
\ref\key We
\by J. Weeks
\paper SnapPea\RM: A computer program for creating and 
studying hyperbolic
\RM3-manifolds
\paperinfo available  by anonymous ftp from {\tt 
geom.umn.edu}
\vol 
\yr 
\pages 
\endref

\endRefs
\enddocument